\documentclass[12pt]{article}
\usepackage{amsmath}
\usepackage{amssymb}
\usepackage{amsfonts}
\usepackage{eucal}
\usepackage[usenames]{color}
\usepackage{graphicx}
\numberwithin{equation}{section}
\newtheorem{Theorem}{Theorem}[section]
\newtheorem{Proposition}[Theorem]{Proposition}
\newtheorem{Definition}[Theorem]{Definition}

\newtheorem{Remark}[Theorem]{Remark}
\begin{document}
\title{The first-crossing area of a diffusion process with jumps over a constant barrier}
\author{Mario Abundo\thanks{Dipartimento di Matematica, Universit\`a  ``Tor Vergata'', via della Ricerca Scientifica, I-00133 Rome, Italy.
E-mail: \tt{abundo@mat.uniroma2.it}}
}
\date{}
\maketitle

\begin{abstract}
\noindent For a given barrier $S$ and a one-dimensional
jump-diffusion process $X(t),$ starting from $x<S,$  we study the
probability distribution of the integral $A_S(x)= \int _0 ^
{\tau_S(x)}X(t) \ dt$ determined by $X(t)$ till its first-crossing
time $\tau_S(x)$ over $S.$ In particular, we show that the Laplace
transform and the moments of $A_S(x)$ are solutions to certain
partial differential-difference equations with outer conditions.
The distribution of the minimum of $X(t)$ in $[0, \tau_S(x)]$ is also
studied. Thus, we extend the results of a previous paper by the
author, concerning the area swept out by $X(t)$ till its
first-passage below zero. Some explicit examples are reported,
regarding diffusions with and without jumps.

\end{abstract}


\section{Introduction}
This paper deals with the first-crossing area, $A_S(x)= \int _0 ^
{\tau_S (x)}X(t) \ dt,$ determined by a one-dimensional
jump-diffusion process $X(t)$ starting from $x,$ till its first-crossing time,  $\tau_S (x) ,$ over a
threshold $S >x,$ and it extends the results of a previous paper by
the author (\cite{abundo:passarea11}), concerning the area swept
out by  $X(t)$ till its first-passage below zero. 
As for results about the integral of $X(t)$ over a deterministic and fixed time interval, see e.g
\cite{abundo:ija08}. \par\noindent
Notice that
$A_S(x)= \overline X _ {[0,\tau_S(x)]} \cdot \tau_S(x),$ where $\overline X _ {[0,\tau_S(x)]}$
denotes the time average of $X(t)$ over the interval $[0, \tau_S(x)].$ \par
The first-crossing area has interesting applications in Biology, for instance in
the framework of diffusion
models for neural activity, if one identifies $X(t)$ with the neuron voltage at time $t,$
and $\tau_S$ with the instant at which the neuron fires, i.e. $X(t)$
exceeds the potential threshold value $S;$ then, $A_S / \tau_S $  represents the time average of the neural
voltage till the first-crossing time over $S.$
Another application can be found in
Queueing Theory, if  $X(t)$ represents the length of a
queue at time $t,$ and one identifies  the first-passage time $\tau_S$ over the threshold $S$
with the overflow time, that is the instant at which the
queue system first collapses;  then, the first-crossing  area represents the
cumulative waiting time experienced by all the ``users''
till the congestion time.
\par\noindent
As for an example from Economics, let us suppose that the
variable $t$ represents the quantity of a commodity that producers
have available for sale and $Y(t)$  describes the price of the
commodity as a function of the quantity in a supply-and-demand
model. Let us assume that $Y(t)$ is
bounded  between $0$ and $S >0,$ and consider the distance of $Y(t)$ from $S,$
that is the process $X(t):=S - Y(t).$  If we denote by $\tau$  the amount of product at
which  $Y(t)$ falls to zero (i.e. $X(t)$ reaches the value $S),$ then
$ \int _0 ^{\tau} X(t) dt = \tau  S - \widetilde A,$ where
$\widetilde A = \int _0 ^\tau Y(s) ds$
provides a measure of the total value that
consumers receive from consuming the amount $\tau$  of the product. \par
In the present article,  we complete the study carried out in \cite{abundo:passarea11};
in fact, for a one-dimensional jump-diffusion process, in place
of the first-passage area below zero, we consider the analogous
problem of first-crossing area over a positive
barrier.\par\noindent Precisely, let be given a barrier $S>0$ and
a one-dimensional jump-diffusion process $X(t),$ starting from
$x<S;$ our aim is  to study the probability distribution of the integral \ \par\noindent
$A_S(x)= \int _0 ^ {\tau_S(x)} X(t) \ dt$ determined by $X(t)$
till its first-crossing time $\tau_S(x)$ over $S.$ We improperly
call $A_S(x)$ ``the first-crossing area of $X(t)$ over $S$''.
Indeed, the area of the plane region determined by the trajectory
of $X(t)$ and the $t$-axis in the first-crossing period $[0,
\tau _S(x)]$ is $\int _0 ^ {\tau _S(x)} |X(t)| \ dt $ which coincides
with $A_S(x)$ only if $X(t)$ is non-negative in the entire interval
$[0, \tau(x)].$
\par Since the topic was studied quite extensively in
\cite{abundo:passarea11} (though in the slight different situation of first passage
below zero), we will omit some details. We will suppose that
$X(t)$ is the solution of a stochastic differential equation of
the form:
\begin{equation} \label{jumpdiffueq}
 dX(t) = b(X(t))dt +\sigma(X(t)) dB_t + \int _{-\infty}^{+\infty} \gamma
(X(t),u) \nu(dt,du)
\end{equation}
with assigned initial condition $X(0) = x  ;$ here $B_t$ is a standard Brownian
motion, $\nu (\cdot ,\cdot)$ is a temporally homogeneous Poisson random
measure (see Section 2 for the definitions), and the functions $b(\cdot), \sigma (\cdot), \gamma (\cdot, \cdot)$
satisfy suitable conditions for the
existence and  uniqueness of the solution (see Section 2).
The coefficients regulate the drift $(b),$ the diffusion $(\sigma),$ and the sizes of the jumps $(\gamma )$ which
occur at (random) exponentially distributed time intervals. \par\noindent
The process $X(t)$ which is the solution of the equation \eqref{jumpdiffueq} reduces to a simple diffusion (i.e. without jumps)
if $\gamma (x, u) = 0,$ and in particular to
Brownian motion with  drift $\mu$ if $ b(x)=  \mu $ and $\sigma (x)=1.$
\par\noindent
Denote by
$$\tau_S(x) = \inf \{ t>0 : X(t) \ge S | X(0)=x \} $$
the first-crossing time over $S$
of the process $X(t)$ starting from $x <S,$ and  assume that $\tau_S(x)$ is finite with probability one.
We will  study the
probability distributions of $\tau_S(x)$  and of the
integral $A_S(x)= \int _0 ^
{\tau_S(x)}X(t) \ dt ,$
as well as their moments; moreover, the distribution of the minimum  of $X(t)$ in $[0, \tau_S(x)]$ will be studied.
In particular, we will show that the Laplace transforms of $A_S(x)$ and $\tau_S(x),$ their moments, as well as the probability
distribution of the
minimum  of $X(t),$ are solutions to certain partial differential-difference equations (PDDEs) with outer conditions.
\par
The paper is organized in the following way: Section 2 contains the statement of the problem and main results, in Section 3 some explicit
examples are reported. Finally, Section 4 is devoted to conclusions and final remarks.

\section{Notations, formulation of the problem and main results}
Let  $X(t) \in I := (\alpha ,\beta ), \ - \infty \le  \alpha < \beta  \le + \infty ,$   a  time-homogeneous, one-dimensional jump-diffusion
process   which  satisfies the
stochastic differential equation (SDE) \eqref{jumpdiffueq}
with assigned initial condition $X(0) = x,$ where $B_t$ is a standard Brownian
motion and $\nu (\cdot ,\cdot)$ is a Poisson random
measure on $(0, \infty) \times \rm I\!R.$ Then, $X(t)$ can be represented as
\begin{equation} \label{integralSDE}
X(t)= x + \int _ 0 ^t b(X(s)) ds + \int _0 ^t \sigma (X(s)) d B_s + \int _0 ^t \int _ { - \infty} ^ {+ \infty } \gamma (X(s), u) \nu (ds,du)
\end{equation}
For the
definitions of the integrals in the right hand side of \eqref{integralSDE} and the Poisson
measure, see  \cite{gimsko:sde72}.
The coefficients $b, \sigma$ and $\gamma$ completely specify
the law of $X(t).$ In particular, an atom $(t,u)$ of the Poisson
random measure $\nu$ causes a jump from $x$ to $x+ \gamma (x,u)$ at time $t,$
if $X(t ^ -)=x.$
We assume that $\nu$ is homogeneous with respect to time translation, that is, its intensity measure $E(\nu(dt, du))$ is of the
form
\begin{equation} \label{poissonmeasure}
E  [ \nu  (dt,du  ) ] = dt  \pi (du)
\end{equation}
for some positive measure $\pi$ defined on $ {\cal B} (\rm I\!R)$ \ $( {\cal B} (\rm I\!R)$ denotes
the Borel $\sigma-$field of subsets of $\rm I\!R),$ and
we suppose that the jump intensity
\begin{equation}
 \Theta = \int _{-\infty}^ {+\infty} \pi (du)  \ge 0
\end{equation}
is finite.
\par\noindent
We make the following assumptions on the coefficients:
\par
\noindent {\bf A1}
$b, \sigma :  I \rightarrow \rm I\!R $ are  continuous
functions  and a constant $K >0$ exists, such that, for
every $ x,y \in I:$
$$ |b(x) - b(y)|  \le K |x-y|$$
$$ b^2(x)+ \sigma ^2 (x)  \le K (1+x^2)$$
\par
\noindent {\bf A2} $ \sigma$ is a non-negative, bounded
function and it is differentiable for every $x$ belonging to the interior of  $I.$
Moreover, there exists a strictly increasing function
$ G : \rm I\!R^+ \rightarrow \rm I\!R $ such that $G (0) =0, \ \int _ {0 ^+}
G ^{-2} (s) ds = + \infty$ and
$$ |\sigma (x) -\sigma (y)| \le G (|x-y|) $$
for every $x,y \in I.$
\par
\noindent {\bf B1} For every $ u$ and $x \in I:$
$$ \int _ {- \infty} ^ {+ \infty} \gamma ^2 (x,u) \pi (du) \le K (1+x^ 2)  $$
\noindent {\bf B2} For every $u$ and $x,y \in I:$
$$ \int _ {- \infty} ^ {+ \infty} | \gamma (x,u) - \gamma (y,u) | \pi (du) \le K |x-y| $$
\par
\noindent
\begin{Remark} {\rm
The conditions  A1,  A2, and  B1, B2 ensure that there exists a
unique non-explosive solution of  \eqref{jumpdiffueq} which is a
temporally homogeneous Markov process (see
\cite{gimsko:sde72}, \cite{ikwa:sde81});   A2 holds, for instance,
if $\sigma (\cdot) $ is Lipschitz-continuous ($G(x)= const \cdot x)$, or
H${\rm\ddot{o}}$lder-continuous of order $ \alpha \in [ 1/2, 1)$ \ ($G(x)= const \cdot x ^ \alpha  ).$  }
\end{Remark}
Notice that, if $\gamma = 0,$  or $\nu = 0,$
then the SDE \eqref{jumpdiffueq} becomes the usual It${\rm \hat{o}}$'s stochastic
differential equation for a simple-diffusion (i.e. without jumps).
\par In the special case when the measure $\pi$ is concentrated
over the set $\{u_1,u_2\}=\{-1,1\}$ with $\pi(u_i)=\theta _i$
and $\gamma(u_i)=\epsilon _i,$ we can rewrite the SDE
\eqref{jumpdiffueq} as
\begin{equation} \label{downuppoisson}
dX(t)=b(X(t))dt + \sigma(X(t)) dB_t + \epsilon _2 dN_2 (t) + \epsilon _1
dN_1 (t)
\end{equation}
where $ \epsilon _1 <0, \ \epsilon _2 >0 $ and $N_ i(t), t \ge 0 $ are independent
homogeneous Poisson processes of
amplitude $1 $ and rates $\theta _1$ and
$\theta _2,$ respectively  governing downward $(N_1)$ and upward $(N_2)$
jumps.

Let $D$ be the class of function $f(t,x)$ defined in $\rm I\!R ^+
\times I,$ differentiable with respect to $t$ and twice
differentiable with respect to $x,$ for which the function
$f(t,x+\gamma(x,u)) - f(t,x)$ is $\pi-$integrable for any $(t,x).$
We recall the generalized It${\rm \hat{o}}$'s formula for
jump-diffusion processes giving the differential of a function $f
\in D$ (see  \cite{gimsko:sde72}):
\begin{equation} \label{itoformula}
df(t,X(t))=\left [ \frac {\partial f}  {\partial t}(t,X(t)) + b(X(t))
\frac { \partial f}  {\partial x} (t,X(t)) + \frac 1  2 \sigma ^2 (X(t))
\frac { \partial ^2 f}  {\partial x ^2 } (t,X(t)) \right ] dt
\end{equation}
$$ + \frac { \partial f } {\partial x}  (t,X(t)) \sigma (X(t)) dB_t
+ \int _ {- \infty} ^ {+ \infty} [f(t,X(t)+\gamma(X(t),u)) - f(t,X(t))] \nu (dt,du) . $$
The differential operator associated to the process $X(t)$ which
is the solution of \eqref{jumpdiffueq},
is defined for any function $f \in D$ by:
\begin{equation} \label{generator}
L f (t,x) = L _d f (t,x) + L _j f (t,x)
\end{equation}
where the ``diffusion part'' is
$$ L _d f (t,x) = \frac 1 2
\sigma ^2 (x) \frac {\partial ^2 f}  {\partial x^2} (t,x)+
b(x) \frac { \partial f}  {\partial x} (t,x)  $$
and the ``jump part'' is
$$L_ j f (t,x) = \int _ {- \infty} ^{+ \infty} [f(t,x+ \gamma(x,u))-f(t,x)] \pi(du) .$$
Then, from \eqref{itoformula}, taking expectation, one obtains
\begin{equation}
E[f(t,X(t))] = f(0,X(0)) + E \left ( \int _0 ^t \left [\frac {\partial f}
{\partial s} (s,X(s)) + L f (s,X(s)) \right ] ds \right ).
\end{equation}
For a barrier $S$ and $x<S,$ we define:
\begin{equation} \label{firstpassagetime}
\tau _S(x) = \inf \{ t >0: X(t) \ge S | X(0)=x \}
\end{equation}
that is the first-crossing time over $S$ of $X^x(t)$ (i.e. the process $X(t)$
starting from $x)$ and  suppose that $\tau _S(x)$ is finite with
probability one. Really, it is possible to show (see  \cite{tuckwell:jap76}, \cite{abundo:pms00}) that the probability $p_0(x)$ that
$X^x(t)$ ever leaves the interval $(- \infty, S)$ satisfies the
partial differential-difference equation (PDDE):
\begin{equation}
L p_0 = 0
\end{equation}
with  outer condition:
$$ p_0(x) =1 \ {\rm if}  \  x \ge S .$$
The equality $p_0(x) =1$ is equivalent to say that $\tau _S(x)$ is finite with probability one.
For diffusion processes without jumps (i.e. $\gamma =0)$
sufficient conditions are also available which ensure that $\tau _S(x)$ is
finite w.p. 1, and they concern the convergence of certain
integral associated to the coefficients of \eqref{jumpdiffueq}
(see Section 3.1 and also  \cite{gimsko:sde72} , \cite{has:sto80}).
\begin{Remark} \label{remark2.2}
{\rm Let us consider the special case when $X(t)= \widetilde X(t) + J(t),$ where $ \widetilde X$ is a simple-diffusion (i.e. without jumps)
and $J$ is a pure-jump process, set $\widetilde m(t) = E(\widetilde X(t)),$ and suppose that
\begin{equation} \label{taufinitojump}
\exists \ \bar t >0 : \widetilde m(t) +E( J(t)) > S, \ {\rm for \ any } \ t \ge \bar t ;
\end{equation}
then, $P( \tau _S(x) < \infty )=1.$
Indeed, $P( \tau _S (x) = \infty ) >0$ implies that, for a set of trajectories having positive probability, it results $X(t) < S$ for
any $t \ge 0.$
Thus, by taking expectation one obtains $E(X(t)) = \widetilde m(t) + E(J(t)) < S$  for any $t \ge 0,$ which contradicts \eqref{taufinitojump}. \par\noindent
If  $J(t) = \epsilon _1 N_1 (t) + \epsilon _2 N_2 (t) $ (cf. \eqref{downuppoisson}), condition \eqref{taufinitojump} becomes
$\widetilde m(t) + (\epsilon_1 \theta _1 + \epsilon_2 \theta _2 )t > S;$ in particular, for $\epsilon_2= \epsilon = - \epsilon _1$
it writes
$\widetilde m(t) > S + \epsilon ( \theta _1 - \theta _2 )t .$ }
\end{Remark}

Let $U$ be a functional of the process $X;$ assume that $\tau _S(x)$ is finite with probability one, and for $\lambda >0$ denote by
\begin{equation}
M_ {U, \lambda} (x) =  E \left [ e^ { - \lambda
\int _0^ {\tau _S(x)} U(X(s)) ds} \right ]
\end{equation}
the Laplace transform of the integral $\int _0 ^ {\tau _S(x)}
U(X(s))ds.$ Then, the following theorem holds (we omit the proof, since it is quite analogous to that of Theorem 2.3 in \cite{abundo:passarea11}):
\begin{Theorem} \label{laplacetransform}
Let $X(t)$ be the solution of the SDE \eqref{jumpdiffueq}, starting from $X(0)=x < S;$
then, under the above assumptions, $M_{ U, \lambda} (x)$ is the solution of
the problem with outer conditions:
\begin{equation} \label{problemlaplace}
\begin{cases}
L M_{ U, \lambda} (x) = \lambda U(x) M_{U, \lambda } (x )  , \  x < S \\
M_ {U, \lambda }(y) =1,  \ {\rm for} \ y \ge S  \\
\lim _ { x \rightarrow - \infty} M_{U, \lambda } (x) =0
\end{cases}
\end{equation}
where $L$ is the generator of $X,$ which is defined by \eqref{generator}.
\end{Theorem}
\par \hfill  $\Box$
\bigskip

We recall that the n-th order moment of $\int _0 ^ {\tau _S(x)}
U(X(s))ds,$ if it exists finite, is given by $(n=1, 2, \dots ):$
$$ T_ n (x) = E \left [  \left ( \int _0 ^ {\tau _S(x)}
U(X(s))ds \right ) ^n \right ]  = (-1)^n \left [ \frac {\partial
^n } {\partial \lambda ^n } M_ {U, \lambda } (x) \right ] _ { \lambda =0} .$$
Then, taking the n-th derivative with respect to $\lambda$ in
both members of the  equation  \eqref{problemlaplace}, and
calculating it for $\lambda =0,$ one easily obtains that
%\begin{Proposition} \label{propositionmoments}
the n-th order moment $T_n(x) \ (n=1, 2, \dots )$ of $\int _0 ^ {\tau _S(x)} U(X(s)) ds, $ whenever it exists finite,  is the solution of the PDDE:
\begin{equation} \label{eqmoments}
LT_n (x) = -n U(x) T_{n-1} (x) , \ \  x < S
\end{equation}
which satisfies
\begin{equation} \label{outerconditionmoments}
 T_n(x)= 0, \ {\rm for } \ x \ge S
\end{equation}
and an  appropriate  additional  condition.

%\end{Proposition}
%\par \hfill  $\Box$

\noindent Indeed, the only condition $T_n (x)=0$ for $x \ge S$ is not sufficient to determinate uniquely the desired solution of
the PDDE \eqref{eqmoments}, because it is a second
order equation. We will return to this problem when we will consider some explicit examples. Note that for a diffusion
without jumps $(\gamma =0 )$ and
for $U(x)\equiv 1,$
 \eqref{eqmoments} is nothing but the celebrated Darling and Siegert's equation  (\cite{darling:ams53}) for the moments of the first-passage time,
 and \eqref{outerconditionmoments}
 becomes simply the boundary condition $T_n(S) =0.$

\begin{Remark} \label{relazconarea}
{\rm In certain cases, the first-crossing time  and the first-crossing area of $X(t)$ over $S$ can be expressed in terms of the first-passage time and the
first-passage area of a suitable  process below zero. For instance, let $X(t) = x + B_t + \mu t ,$ with $\mu >0,$
i.e. Brownian motion with positive drift $\mu;$ then
$\tau _S(x) = \inf \{ t>0 : X(t) \ge S  \}$ has the same distribution as
$\inf \{ t >0 : S-x + B_t - \mu t \le 0 \} \equiv \widetilde \tau (S-x),$ where $\widetilde \tau (y)$ denotes the
first-passage time of the process $Y(t)=y +B_t - \mu t $ below zero. Moreover, it is easy to see that
\begin{equation} \label{areafunctionareatilde}
A_S(x) = \int _0 ^{ \tau_S (x) } X(t) dt = S \cdot  \widetilde \tau (S-x) - \widetilde A (S-x)
\end{equation}
where $\widetilde A (y)$ denotes the
area swept out by $Y(t)$ till its
first-passage below zero. \par\noindent
Thus, the  moments of $\tau _S(x)$ can be soon obtained by those of $\widetilde \tau (S-x),$ while
by using \eqref{areafunctionareatilde}, the mean of $A_S(x)$ can be obtained in terms of
 $E( \widetilde \tau (S-x))$ and $E( \widetilde A (S-x))$ (the first two moments of $\widetilde \tau$ and $ \widetilde A$ were obtained in \cite{abundo:passarea11} ). Notice that,
the computation of $E(A_S(x)^2)$ requires instead also the knowledge of the covariance  of $\widetilde \tau $ and $ \widetilde A .$
}
\end{Remark}
\bigskip

\noindent {\bf Distribution of the minimum of $X(t).$} \par\noindent
Now, we will study the probability distribution
of the maximum downward displacement (i.e. the minimum) of the jump-diffusion $X(t)$ starting
from $x <S,$  till its first passage over $S,$ that is,
$\mathfrak{m}(x) = \min _ { t \in [0, \tau _S(x)] } \{ X(t) | X(0)=x \}.$ For any $z \le x \le S$
the event $\{ \mathfrak{m}(x) > z \}$ is nothing but the event ``$X(t)$
first exit the interval $[z,S]$ through the right end $S$''; so, by
the well-known result about the exit probability of a
jump-diffusion from the right end of an interval (see
\cite{abundo:pms00}), we obtain that $v(x):= P(\mathfrak{m}(x) >z),$ as a function of $x,$ is solution of the
equation
$ Lv =0$ with conditions $  v(y) =1, \ y \ge S; \ v(y)=0, \ y \le z .$
Thus, we get:
\begin{Proposition} \label{thmaxdisplacement}
$P( \mathfrak{m}(x) \le  z )$ is the solution of the problem with outer conditions:
\begin{equation} \label{maxdisplacement}
\begin{cases}
L w (x) = 0  , \ x \in (z,S) \\ w(y)= 0 , \ y \ge S  \\ w(y)= 1, \ y \le z
\end{cases}
\end{equation}
\end{Proposition}
\par \hfill  $\Box$

\section{A few examples}
In this section we will compute explicitly the moments of $\tau_S(x)$ and  those of the first-crossing area $A_S(x)$ for certain
jump-diffusion processes. We start with considering diffusions without jumps.
\subsection{Simple diffusions (i.e with no jump)}
Let $X(t)$ be the solution of \eqref{jumpdiffueq}, with $\gamma \equiv 0,$ that is:
\begin{equation}
\begin{cases}
dX(t) = b(X(t))dt +\sigma(X(t)) dB_t \\ X(0) =x
\end{cases}
\end{equation}

\noindent In this case $\tau _S(x) \equiv \inf \{ t >0: X(t)=S | X(0)=x \} ,$ that is the first-passage time of $X^x(t)$ through $S.$ \par\noindent
Let us consider the functions ($c$ is a constant) :
\begin{equation} \label{phieq}
 \phi(x) = \exp \left ( - \int _c ^x \frac {2b(s) } {\sigma^2 (s) } ds \right )
\end{equation}
\begin{equation} \label{xieq}
 \xi(x) = \phi (x) \int _c ^x \frac {2} { \sigma^2 (s) \phi(s) } ds .
\end{equation}
As it is well-known (see e.g \cite{gimsko:sde72}),  a sufficient condition in order that $\tau _S(x)$ is finite with probability one,
namely the boundary $S$ is attainable, is that the function $\xi(x)$ is
integrable in a neighbor of $S.$\par
Since the generator $L$ coincides with its diffusion part $L_d,$
by Theorem \ref{laplacetransform} we obtain that, for $x <S,$
$M_{ U, \lambda} (x) =  E \left [ e^ { - \lambda
\int _0^ {\tau _S(x)} U(X(s)) ds} \right ]$ is the solution of the problem with boundary conditions ($M'$ and $M''$ denote first and second derivative with
respect to $x)$ :
\begin{equation} \label{laplaceproblemsimple}
\begin{cases}
\frac 1 2 \sigma ^2 (x) M''_ {U, \lambda } (x) + b(x) M'_ {U, \lambda } (x) =
\lambda U(x) M_ {U, \lambda } (x ) \\
M_ {U, \lambda }(S) =1 \\
\lim _ {x \rightarrow - \infty} M_ {U, \lambda } (x) =0
\end{cases} .
\end{equation}
Moreover, by \eqref{eqmoments}, \eqref{outerconditionmoments}  the n-th order moments $T_n(x)$ of $ \int _0 ^ {\tau(x)}
U(X(s))ds,$ if they exist,
satisfy the recursive ODEs:
\begin{equation} \label{eqmomentsimple}
\frac 1 2 \sigma^2(x) T''_n(x) + b(x) T'_n(x) = -n U(x) T_{n-1} (x), \ {\rm for} \  x < S
\end{equation}
with the condition $T_n(S)= 0,$ plus an  appropriate  additional  condition.
\par
\noindent Finally, as regards the minimum  $\mathfrak{m}(x),$ it turns out that its  distribution
$F(z)= P(\mathfrak{m}(x) \le z)$ is
the solution of the problem with boundary conditions:
\begin{equation} \label{maxdisplacementsimple}
\begin{cases}
L w (x) = 0  , \ x \in (z,S) \\ w(z) =1  \\ w(S) =0
\end{cases} .
\end{equation}
\newpage
\noindent {\bf Example 1} (Brownian motion with drift $ \mu >0)$ \par
\noindent Let be $X(t)= x + \mu t + \sigma B_t ,$ with $\mu , \sigma >0$ and $0< x < S.$
Without loss of generality, we can assume $\sigma =1$  (otherwise, dividing by $\sigma$ one reduces to  this case).
Note that, since the drift is positive, $\tau_S(x)= \tau_S^ \mu (x) = \inf \{t>0 | X(t)=S \} $ is finite with probability one,
for any $x < S.$
Taking $b(x) =  \mu, \ \sigma(x)=1,$ the equation in \eqref{laplaceproblemsimple}
for $M_ {U, \lambda } (x) = E(e^ {- \lambda \int _0 ^ { \tau _S ^ \mu (x) }
U(X(s) ds } )$
becomes
\begin{equation} \label{brownianlaplace}
\frac 1 2 M'' _ {U, \lambda } (x) + \mu M'_ {U, \lambda } (x) -  \lambda U(x) M_{ U, \lambda} (x) =0 .
\end{equation}
{\bf (i)} The moment generating function  of $\tau_S ^ \mu(x)$
\par \noindent By solving \eqref{brownianlaplace}  with $U(x)=1,$ we
explicitly obtain:
$$ M_ {U, \lambda } (x) = c_1 e^ { \rho _1 x} + c_2 e ^ {\rho _2 x }$$
where $ \rho _1 = - \mu - \sqrt { \mu ^2 + 2 \lambda } < 0, \ \rho _2 = - \mu + \sqrt { \mu ^2 + 2 \lambda } >0;$
the constants $c_1$ and $c_2$ must be determined by the boundary conditions. Indeed,
$M_{ U, \lambda} ( - \infty) =0$ gives $c_1 =0,$ while
$M_ {U, \lambda } (S)=1$ implies $c_2= e^{ (\mu - \sqrt { \mu ^2 + 2 \lambda })S}.$
Thus, we get:
\begin{equation}
M_ {U, \lambda } (x)= E(\exp( -\lambda \tau_S ^ \mu (x))= \exp [( \mu - \sqrt { \mu ^2 + 2 \lambda } )(S-x)] .
\end{equation}
This Laplace transform can be explicitly inverted (see  \cite{kartay:sto75}), so obtaining the well-known expression of the density
of $\tau_S ^\mu(x):$
\begin{equation} \label{brownianfirstpassagedensity}
f_ {\tau _S ^\mu(x)} (t) = \frac {S-x } {\sqrt { 2 \pi } t^ { 3/2} } e^ {-(S-x - \mu t)^2/ 2t } .
\end{equation}
For $\mu >0$ the moments $T_n (x)= E(\tau_S^ \mu(x))^n$ of any order $n,$ are finite and they can be easily obtained by calculating
$ (-1)^n  [\partial ^n  M _ {U, \lambda } (x) / \partial \lambda ^n]_ { \lambda =0}.$
We obtain, for instance:
\begin{equation} \label{momentstauforBM}
E(\tau_S ^\mu(x))= \frac {S- x} \mu , \ \
E( (\tau _S^\mu (x))^2)= \frac {S-x} { \mu ^3} + \frac {(S-x)^2} { \mu ^2}
\end{equation}
 (cf. Remark \ref{relazconarea}
and the results for BM with negative drift
in \cite{abundo:passarea11} ).
Note that $E(\tau _S ^\mu(x)) \rightarrow + \infty,$ as $ \mu \rightarrow 0.$
As easily seen, for any $x <S$ it results $M_ {U, \lambda } (x) \rightarrow 1,$ as $\mu \rightarrow + \infty,$ or, equivalently,
$\tau _S ^\mu(x)$ converges to $0$ in distribution, and so $E((\tau _S ^\mu(x))^n) \rightarrow 0,$ as $\mu \rightarrow + \infty,$
for any $n$ and $x <S.$
\bigskip

\noindent {\bf (ii)} The moments of $A _S^\mu(x)= \int _0 ^ { \tau _S ^ \mu (x)} (x+ \mu t + B_t) dt $ \par
\noindent For $U(x)=x$ the equation \eqref{brownianlaplace} becomes
\begin{equation} \label{eqlaplacearea}
\frac 1 2 M''_ {U, \lambda } + \mu M' _ {U, \lambda } = \lambda x M_ {U, \lambda }
\end{equation}
with conditions $M_{U, \lambda } (S)=1, \ M_{U, \lambda } (- \infty ) =0;$
now $M_{ U, \lambda} (x) = E( e ^ { - \lambda A _S^ \mu (x) } ) .$
%If $\mu$ is replaced with $- \mu$ one obtains the Schrodinger equation for a quantum particle moving in a
%uniform field (see e.g \cite{kearney:jph05},
%\cite{kearney:jph07}).
Unfortunately, its explicit solution
cannot be found in terms of elementary functions, but it can be
written in terms of the Airy function (see \cite{kearney:jph05},
\cite{grand:tab80}) though it is impossible to
invert the Laplace transform $M_ {U, \lambda }$  to obtain the
probability density of $A _S ^ \mu (x).$ \par\noindent In the special case
$\mu =0,$ it can be shown (see \cite{kearney:jph05}, \cite{abundo:passarea11}) that the solution
of \eqref{eqlaplacearea} is:
\begin{equation} \label{laplacetransformarea}
M_{ U, \lambda} (x) = 3 ^ {2/3} \Gamma \left ( \frac 2 3 \right ) {\rm
Ai} (2 ^ {1/3} \lambda ^ {1/3} (S-x) )
\end{equation}
where ${\rm Ai }(z)= \pi ^{-1} \sqrt { z/3} \  K_{1/3} \left ( \frac 2 3 z ^{3/2} \right )$ denotes the Airy function, and
 $K _ \nu (z)$ is a modified Bessel function (see \cite{abramo:handbook65}).
Calculating the derivative with respect to $\lambda $ in \eqref{laplacetransformarea}, we obtain:
$$  \frac \partial { \partial \lambda } M_{U, \lambda } (x) = \left ( \frac 2 3 \right ) ^ {1/3} \Gamma \left ( \frac 2 3 \right )  (S-x)
 {\rm Ai } '(2 ^ {1/3} \lambda ^ {1/3} (S-x) ) \cdot \lambda ^{ - 2/3}  $$
By using the fact that ${\rm Ai}  '(0) = - \left ( 3 ^ {1/3} \Gamma \left ( \frac 1 3 \right ) \right ) ^{-1} $ (see \cite{abramo:handbook65}), it
follows that $\frac \partial { \partial \lambda } M_{U, \lambda } (x) | _{\lambda =0 } = - \infty ,$  namely, the expectation of
$A _S^0(x)$ is infinite. Notice that, unlike the Laplace transform inversion reported in \cite{abundo:passarea11}
for a similar case (see equation (3.12) therein), we are not able to invert explicitly the Laplace transform \eqref{laplacetransformarea} to find the density of
$A_S^0;$ in fact, in the case considered in \cite{abundo:passarea11} this was possible, thanks to an integral identity
(see e.g. \cite{kearney:jph05}, \cite{grand:tab80} ) involving the modified Bessel function, which serves the purpose when the support of the
candidate density of $A_S^0 $ is $(0, + \infty ).$
%; then, by  inverting
%this Laplace transform one finds that the first-crossing area
%density is (\cite{kearney:jph05}):
%\begin{equation}
%f_ {A^ 0 _S (x)} (a) = \frac {2^ {1/3} } {3 ^ {2/3} \Gamma ( \frac 1 3 )
%} \ \frac {x } {a ^ {4/3} } \ e^ { -2 x^3 / 9a } \ .
%\end{equation}
%Thus,  the distribution of $A _S^ 0(x)$ has an algebraic
%tail of order $ \frac 4 3 $ and so the moments of all orders are
%infinite. \par
\par
For $\mu >0,$
we will find closed form expression for the first two moments of $A _S ^ \mu(x),$ by solving \eqref{eqmomentsimple} with $U(x) =x.$
For $n=1,$ we get that $T_1 (x)= E(A _S ^\mu(x))$ must satisfy the equation:
\begin{equation} \label{meanarea}
\begin{cases}
 \frac 1 2 T_1''(x) +  \mu T_1'(x) = - x \\ T_1(S)=0
\end{cases}
\end{equation}

Of course, the only  condition $T_1(S)=0$ is not sufficient to uniquely determinate the solution.
The general solution of \eqref{meanarea} involves  arbitrary constants $c_1$ and $c_2$ and, as easily seen, it is given by
$$T_1(x)= c_1 + c_2 e ^ {-2 \mu x } - \frac {x^2 } {2 \mu } + \frac {x } {2 \mu ^2 }  \ .$$
By imposing  that $T_1(S)=0$ and that for any $x \le 0$ it must be $ T_1(x) \rightarrow 0,$ as $\mu \rightarrow + \infty,$
we find that the mean first-crossing area is
\begin{equation} \label{expressionmeanarea}
E(A _S ^ \mu (x))= \frac { S-x} {2 \mu } \left [S+x - \frac 1 \mu \right ]
\end{equation}
This formula also follows by taking expectation in \eqref{areafunctionareatilde}, and using that
$E( \widetilde \tau (y)) = \frac y  \mu ,$ \par\noindent
$E( \widetilde A (y)) =
\frac { y^2 }  {2 \mu} + \frac y  { 2 \mu ^2} $ (see \cite{abundo:passarea11}).\par
As far as the second moment of $A _S ^ \mu (x)$ is concerned, we have to solve \eqref{eqmomentsimple} with $U(x)=x$ and  $n=2,$
obtaining the equation for $T_2(x)= E [ ( A _S^ \mu (x)) ^2 ]:$
\begin{equation} \label{secondarea}
\begin{cases}
 \frac 1 2 T_2''(x) +  \mu T_2'(x) = - \frac x  \mu (S-x) (S+x - 1/ \mu ) \\ T_2(S)=0
\end{cases}
\end{equation}
As before, the only  condition $T_2(S)=0$ is not sufficient to uniquely determinate the solution.
The general solution of \eqref{secondarea}, which involves two arbitrary constants $c_1$ and $c_2,$  is given by
$$T_2(x)= c_1 + c_2 e ^ {-2 \mu x } + A x^4 +B x^3 + C x^2 + Dx ,$$
where $$A = \frac 1 {4 \mu ^2 } , \ B = - \frac 5 { 6 \mu ^3 }, \ C = \frac S {4 \mu ^4 } (2 \mu +1 - 2 \mu ^2 S ), \
D = - \frac S {4 \mu ^5 } (2 \mu +1 - 2 \mu ^2 S ).$$
By imposing that, for any $x \le 0$ it must be $T_2(x) \rightarrow 0,$ as $ \mu \rightarrow + \infty ,$ we find $c_2 =0;$
moreover, by $T_2(S)=0$ we get $c_1=  - (AS^4 + B S^3 + C S^2 +DS) .$
Thus the second order moment of the first-crossing area is
\begin{equation} \label{expressionsecondarea}
E[(A _S^\mu (x))^2]= A (x^4 - S^4) + B (x^3 - S^3 ) + C (x^2 - S^2 ) + D(x-S),
\end{equation}
where the constants $A, B, C, D$ are as above.
Finally, by \eqref{expressionmeanarea} and \eqref{expressionsecondarea}, one can obtain the variance of $A _S^ \mu (x).$
\par\noindent
Notice that $E[(A_S ^\mu (x)) ^2 ]$ has to be the only non-negative solution to
\eqref{secondarea}.  Indeed, expression \eqref{expressionsecondarea} loses meaning if it becomes negative
for some choice of $x, \ \mu$ and $S;$ in that case, the second moment of $A_S ^\mu (x)$ does not exist.
\par\noindent
Since a closed form expression for the density of the first-crossing area $A _S^ \mu (x)$
cannot be found for $\mu >0 ,$ it must be obtained numerically. As in the case of the first-passage area below zero
(see  \cite{abundo:passarea11}, \cite{abundo:browarea11}),
we  have estimated it  by simulating a large number of trajectories of Brownian motion
with drift $ \mu >0,$ starting from the initial state $x>0.$
The first and second order moments of the first-crossing time $\tau_S ^\mu $ and
of the first-crossing area $A_S^\mu $ thus obtained, well agree with the exact
values. For $S=2, \ x =1$ and several values of  $\mu >0,$ we
report in the Figure 1 the estimated density of the
first-crossing area, regarding $X(t)= x + \mu t + B_t.$

\begin{figure}
\centering
\includegraphics[height=.4\textheight]{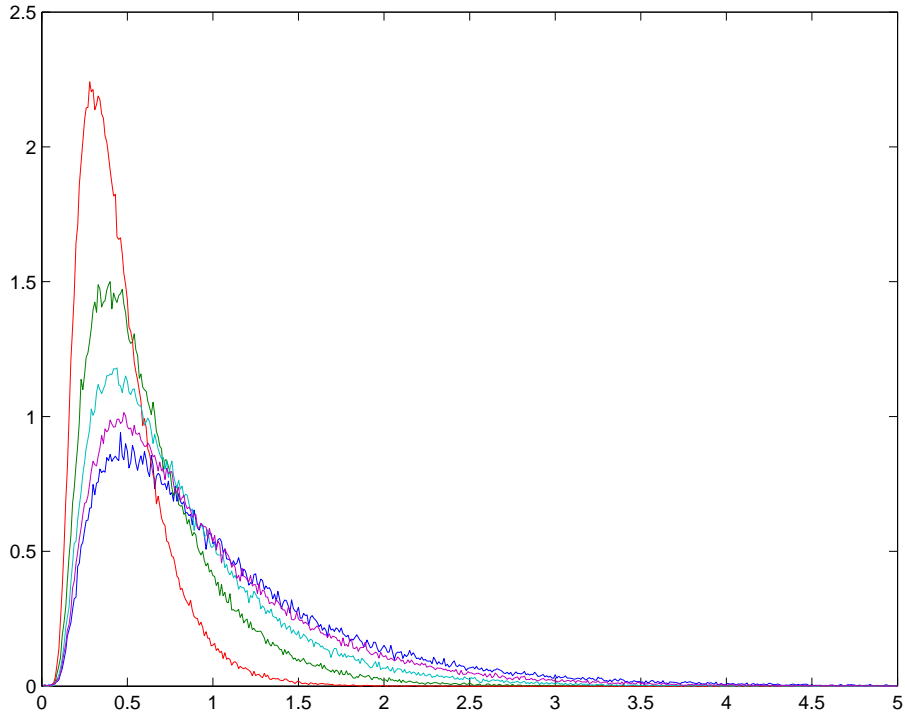}
\caption{Estimated density of the first-passage area
$A_S^ \mu (x)$ of $X(t)= x + \mu t + B_t$ over $S,$ for $S=2, \ x=1$ and several values of $\mu.$ From  top to bottom,
with respect to the peak of the curve: $\mu = 3; \ \mu = 2;  \ \mu
= 1.5;  \  \mu = 1.2 ;  \ \mu = 1.$
 }
\end{figure}

\noindent For some values of parameters we have compared
the estimated density of the first-crossing area with a suitable
Gamma density; in the Figure 2, we report for $S=2, \ x=1$ and $\mu = 1.5,$
the comparison of the Laplace transform $M_ \lambda (x)= E(e^ {-
\lambda A _S ^\mu (x)} )$ of the estimated first-crossing area density  and the
Laplace transform of the Gamma density with the same mean and
variance. Although the two curves agree very well for small values of $\lambda  >0$
(this implying a good agreement between the
moments), for large values of $\lambda $ the graph of the Laplace
transform of the estimated density of $A_S^\mu (x)$ lyes below the other one, which is compatible
with an algebraic tail for the distribution of $A_S^\mu (x).$

\begin{figure}
\centering
\includegraphics[height=.4\textheight]{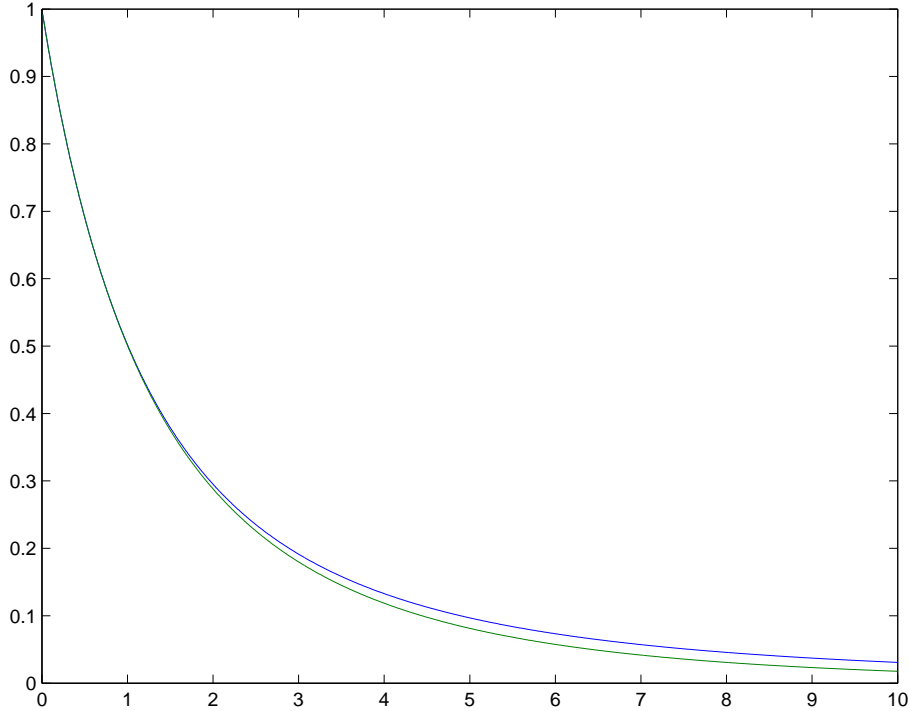}
\caption{Comparison of the Laplace transform $M_ \lambda (x)= E(e^ {- \lambda A_S ^\mu(x)} )$ of the estimated density of the first-passage
area $A_S ^\mu(x)$ (lower curve) and the Laplace transform of the Gamma density with the same mean and variance (upper curve),
as functions of
$\lambda,$  for $S=2, \ x=1$ and $\mu = 1.5.$
}
\end{figure}
\bigskip

\noindent {\bf (iii)} The distribution  of the minimum  $\mathfrak{m} ^\mu (x).$\par
\noindent Set $$\mathfrak{m} ^\mu (x)= \min _{[0,  \tau _S ^\mu (x) ] } \{ x + \mu t + B_t \}.$$ Its distribution function $P(\mathfrak{m} ^\mu(x) \le z )= w(x)$ is the solution of the problem with boundary conditions:
\begin{equation}
\begin{cases}
\frac 1 2 w''(x) + \mu w'(x) =0 \\ w(z)=1, \ w(S)=0
\end{cases}
\end{equation}
By solving the above equation, we obtain, for $z \le x:$
\begin{equation} \label{dstributionminimum}
P(\mathfrak{m} ^\mu (x) \le z ) = w(x) = \frac {e^ { -2 \mu x} - e ^ { 2 \mu S } } {e^{- 2 \mu z} - e^{-2 \mu S } }
\end{equation}
\noindent Then, calculating the derivative with respect to $z,$ we get the probability density of $\mathfrak{m} ^\mu (x):$
\begin{equation}
 f_ {\mathfrak{m} ^\mu (x)} (z)  = \frac {2 \mu e ^ { -2 \mu z } \left ( e^ { -2 \mu x } -e^{-2 \mu S} \right ) }
{ \left ( e^{-2 \mu z } - e^{-2 \mu S } \right ) ^2 } \ {\bf 1} _ { ( - \infty , x ]} (z)
\end{equation}
For $\mu =0,$ taking the limit in the above expression, we get $f_{\mathfrak{m} ^0 (x)} (z) = \frac {S-x} {(S- z)^2 } \ {\bf 1} _ { ( - \infty , x ]} (z) ,$ and so
the moments of $\mathfrak{m}^0 (x)$ of all
orders $n\ge 1$ are
infinite.\par
\noindent On the contrary, for $\mu >0 $ the minimum  $\mathfrak{m} ^\mu (x)$ possesses finite moments of all orders,
because $f_ {\mathfrak{m} ^\mu (x)} (z)$ behaves like $const / e^{ -2 \mu z }, $ for $z \rightarrow - \infty .$

\begin{Definition}
We say that a one-dimensional diffusion $X(t)$ with $X(0)=x,$ is conjugated to BM if there exists an increasing differentiable
function $u(x),$ with $u(0)=0,$
such that $X(t)=u^{-1} \left ( B_t + u(x) \right ).$
\end{Definition}

\begin{Remark} {\rm
\noindent If $X(t)$ is conjugated to Brownian motion via the function $u,$ then for $x < S:$
$$ \tau _S(x)= \inf \{ t >0: X(t)=S | X(0)=x \} = \inf \{ t >0 : u(X(t)) =u(S) |u(X(0))= u(x) \} $$
$$ = \inf \{ t >0: B_t + u(x) =u(S) \} = \tau_{S'}^B (u(x)) $$
where $S'= u(S)$ and $\tau_{S'} ^B (y)$ is the first hitting time to $S'$ of $y + B_t$ (i.e. BM without drift starting from $y).$
Thus, the first-passage time through $S$ of the process $X(t)$ starting from $x <S,$ is nothing but the first hitting time to $S'$ of
BM starting from $u(x),$
and $$A_S(x)= \int _0 ^ { \tau _S(x)} X(t) dt = \int _0 ^ { \tau_{S'}^B (u(x)) } u^{-1} (B_t + u(x)) dt .$$
Notice that, though $\tau _S(x)$ turns out to be finite with probability one,
it results $E(\tau _S(x))=  + \infty .$ Moreover
$$ M_ {U, \lambda } (x)= E \left [ \exp \left ( - \lambda \int _0 ^ {\tau _S(x)} U(X(t)) dt \right )  \right ] =
E\left [ \exp \left ( - \lambda \int _0 ^ {\tau _{S'}^B(u(x))} V(B_t) dt \right ) \right ) ,$$
where  $V(s)= U \left (u^{-1}  ( s + u(x)  ) \right ).$ \par\noindent
Therefore, the Laplace transform $M_  {U, \lambda } (x)$ of $\int _0 ^ {\tau _S(x)} U(X(t)) dt,$ associated to the functional $U$ of the process $X,$
is nothing but the Laplace transform
$M_ {V, \lambda}  ^{S'} (y)$ of $\int _0 ^ {\tau_{S'} ^B
(y)} V(B_t) dt,$  associated to the functional $V$ of BM, where $y=u(x).$ \par\noindent
This means that the equation \eqref{laplaceproblemsimple} is easily reduced to the analogous equation for BM starting from $u(x),$ with $U$ replaced by $V.$
 }
\end{Remark}

\noindent {\bf Example 2.} \par\noindent
A class of diffusions  conjugated to BM is given by processes $X(t)$ which are solutions of SDEs such as
\begin{equation}
 dX(t) = \frac 1 2 \sigma (X(t)) \sigma ' (X(t)) dt + \sigma (X(t)) dB_t, \  X(0)=x_0 \label{conjdiffu}
 \end{equation}
with $\sigma (\cdot) \ge 0.$
Indeed, if the integral
$v(x) \doteq  \int _{x_0} ^ x \frac {1}  {\sigma (r) }dr $
is convergent for every $x ,$
by It${\rm \hat o}$'s  formula, we obtain that $X(t)= v^ {-1} (B_t +v(x_0)).$
\bigskip

\noindent {\bf (i)}  \ (Feller process or Cox-Ingersoll-Ross (CIR) model) \par\noindent
For $a ,b \ge 0, $  let $X(t)$ be the solution of the SDE
\begin{equation} \label{fellereq}
dX(t) = (a+ b X(t)) dt + \sqrt {X(t) \vee 0} \ dB_t \ , X(0) = x_0 \ge 0
\end{equation}
(note that, although $\sqrt x$ is not Lipschitz-continuous, the
solution is unique because $\sqrt x$ is H${\rm
\ddot{o}}$lder-continuous of order $\frac 1 2$ (see condition
A2)). The process $X(t)$ turns out to be  non-negative for all
$t\ge 0$ (see  \cite{abundo:saa06}, \cite{abundo:pms97}). If $b=0$ and $ a = \frac
1 4,$ the SDE \eqref{fellereq} becomes:
$$ dX(t) = \frac 1 4  dt + \sqrt {X(t) \vee 0} \ dB_t \ , X(0)=x_0 $$
and $X(t)$ turns out to be conjugate to BM via the function $v(x) =2 \sqrt x$ i.e. $X(t) = \frac 1 4 (B_t +2 \sqrt { x_0} )^2;$ the SDE is
obtained by taking $\sigma (x)= \sqrt { x \vee 0 }$ in \eqref{conjdiffu}. \bigskip

\noindent {\bf (ii)} (Wright \& Fisher-like process).\par\noindent
The diffusion described by the
SDE:
\begin{equation} \label{wrightfishereq}
dX(t) = (a + b X(t)) dt + \sqrt { X(t) (1-X(t))\vee 0   } \ dB_t, \  X(0)=x_0 \in [0,1]
\end{equation}
with $a \ge 0$ and $a + b \le 0$ does not exit from the
interval $[0,1]$ for any time (see \cite {abundo:saa06}, \cite{abundo:pms97}). This equation is used for instance
in the Wright-Fisher model for population genetics and in certain
diffusion models for neural activity \cite{lansky:jtb94}.
For $a = \frac 1 4$ and $b = - \frac 1 2,$
$X(t)$ turns out to be conjugated to BM
via the function
$v(x) = 2 \arcsin \sqrt {x},$ i.e. $X(t) = \sin ^2 ( B_t /2 + \arcsin \sqrt {x_0}).$
\par\noindent
For these special values of parameters, the SDE \eqref{wrightfishereq} becomes:
$$dX(t) = \left ( \frac 1 4 - \frac 1 2 X(t) \right ) dt + \sqrt {X(t)(1-X(t))} \ dB_t \ , X(0) = x_0 \in [0,1] $$
and it is obtained from \eqref{conjdiffu} by taking $\sigma(x)= \sqrt { x(1-x)\vee 0 } \ .$
Notice that, for $0 \le x < S \le 1,$   it results $A_S(x)= \int _0 ^ {\tau _S(x)} X(t) dt  \le \tau_S(x), $
since $0 \le X(t) \le 1.$\bigskip

\noindent {\bf Example 3} (Ornstein-Uhlenbeck process) \par
\noindent Let  $X(t)$ be
the solution of the SDE:
\begin{equation} \label{ouequation}
dX(t) = - \mu X (t) dt + \sigma dB_t, \ X(0)= x
\end{equation}
where $\mu$ and $ \sigma$ are positive constants.
By calculating the functions $\phi$ and $\xi$ in \eqref{phieq} and \eqref{xieq}, we obtain:
$$ \phi (x)= const \cdot \exp \left ( \frac { \mu x^2} {\sigma ^2 } \right ) $$
and
$$ \xi (x) = const \cdot \frac 2 \sigma \sqrt { \frac \pi \mu } \exp \left ( \frac { \mu x^2} {\sigma ^2 } \right ) \left [ \Phi \left ( \frac {x \sqrt { 2 \mu } } {\sigma } \right )
- \frac 1 2 \right ] $$
where $ \Phi$ denotes the distribution function of the standard Gaussian variable.
Since $\xi (x)$ is integrable in a neighbor of $x=S,$ the boundary $S$ is attainable. \par
\noindent The explicit
solution of \eqref{ouequation} is
$ X(t) = e ^ { - \mu t} \left ( x + \int _0 ^t \sigma e ^ {\mu s} dB_s \right )$ (see e.g. \cite{abundo:asy09}, \cite{abundo:euc07}).
By using a
time--change, we can write
$ \int _0 ^t \sigma e ^ {\mu s} dB_s =  B _ {\rho (t)},$ where
$ \rho (t) = \frac { \sigma ^2 }  { 2 \mu} \left ( e ^ {2 \mu t} -1 \right ).$
Then, one gets $X(t) = e^ { - \mu t } (x + B ( \rho (t)).$  Let us consider e.g. the time dependent boundary
$S(t)= H(t)= \alpha e^{ - \mu t },$ then:
$$ \tau_H(x)= \inf \{ t >0 : X(t)=H(t) | X(0)=x \} = \inf \{ t >0 : B(\rho(t)) + x  =\alpha \} $$
and so
$$ \rho ( \tau _H (x)) = \inf \{ s >0 : x  + B_s =\alpha \} = \tau _ \alpha ^B (x),$$
that is $\tau _H(x)= \rho ^ {-1}  ( \tau _ \alpha  ^B ( x ) ),$ where $\tau _ \alpha  ^B (x)$ is the first-passage time through
$\alpha$  of Brownian motion starting from $x.$
Thus, $\tau _H (x)$ has density
$$ f_ { \tau _H (x)} (t) = f _ {\tau _ \alpha ^B (x)}  ( \rho (t)) \rho '(t) $$
where $f _ {\tau _ \alpha ^B (x)}$ denotes the density of $ \tau _ \alpha  ^B ( x  ),$ which is
given by \eqref{brownianfirstpassagedensity} with $\mu =0$ and $ S = \alpha .$
Notice that even the
mean of $\tau _H (x)$ is impractical to be directly calculated  by using $f_{\tau _H (x)} (t).$ \par\noindent
However, for a constant boundary $S$ it results (see e.g. \cite{ricc:jap99}):
$$ E( \tau_S(x))= \frac 1 \mu \big \{ \sqrt \pi \left [ \phi_1 \left ( S \sqrt \mu  / \sigma \right )-
\phi_1 \left ( x \sqrt \mu / \sigma  \right ) \right ] + \psi_1 \left ( S \sqrt \mu / \sigma \right )-
\psi_1 \left ( x \sqrt \mu / \sigma \right )  \big \} ,$$
where
$$ \phi _1(z)= \int _0 ^z e^{t^2} \ dt = \sum _ {k=0}^\infty \frac {z^{2k+1} } {(2k+1)k! } , $$
$$ \psi_1 (z)= 2 \int _0 ^z \ du \ e^ {u^2} \int _0 ^u \ dv \
e^ {- v^2 } = \sum _ {k=0} ^\infty \frac {2^k } {(k+1)(2k+1)!! } z^{2k+2} . $$
%  !!!!CHECKED !!!!
Notice that, if $S>0,$ as seen by using Hospital's rule, 
$ \ E(\tau _S(x))$ tends to infinity, for $ \mu \rightarrow
0,$ since the OU process
reduces to BM; also for $ \mu \rightarrow + \infty $ it tends to infinity,
which is indeed natural, since the drift tends to $ - \infty . $
\par\noindent As far as the minimum
 $\mathfrak{m}(x)$ is concerned, by Proposition
\ref{thmaxdisplacement} its distribution $w(x)=P( \mathfrak{m}(x) \le z)$
satisfies the problem:
$$ \begin{cases}
\frac 1 2  \sigma ^2 w''(x) -  \mu x w'(x) =0 , \ z < x < S \\
w(z)=1, \ w(S) =0
\end{cases}
$$
whose solution is:
$$ P( \mathfrak{m}(x) \le z )= \left ( \int _x^S e^ { \frac { \mu } {\sigma ^2 } t^2 } dt \right )
\left (\int _z^S e^ { \frac { \mu} {\sigma ^2 } t^2 } dt \right ) ^ {-1} , \ z \le x  .$$
Thus, the density of $\mathfrak{m}(x)$ is:
$$ \frac {d} {dz} P( \mathfrak{m}(x)) \le z )= \left [
e ^ { \frac { \mu} {\sigma ^2 } z^2 } \int _ x ^S  e ^ { \frac { \mu} {\sigma ^2 } t^2 } dt
\left (  \int _ z ^S  e ^ { \frac { \mu} {\sigma ^2 } t^2 } dt \right )^ {-2} \right ] \cdot {\bf 1} _ { (- \infty, x ]} (z).
$$

\subsection{Diffusions with jumps}
{\bf Example 4} (Poisson process) \par
\noindent
For $x>0,$ let us consider the jump-process $X(t)= x + N_t,$ where $N_t$ is a homogeneous Poisson process with intensity $\theta >0,$ namely $N_0=0$
and its
jumps, of amplitude $1,$
occur at independent instants, exponentially distributed with parameter $\theta.$
This means that
$$P( N_t = k) = \frac { e^{ - \theta t} ( \theta t ) ^k}  { k!}, \ {\rm for} \ k=0, 1, \dots  $$
The infinitesimal generator of the process is defined by:
$$L g(x) = \theta [g(x+1) - g(x)] , \ g \in C^0(\rm I\!R) $$ and, for $S>0$ and $x < S,$
$ \tau _S (x) = \inf \{ t >0 : x + N_t \ge  S \}.$ \par

By Theorem \ref{laplacetransform} with $U(x)=1,$ it follows that the
Laplace transform $M_ {U, \lambda } (x)$ of $\tau _S(x)$ is the solution of the equation $L M _ {U,  \lambda } (x) = \lambda M _ {U,  \lambda } (x),$ with
outer condition $M_ {U, \lambda } ( y) =1$ for $ y \ge S.$ By solving this equation, we get:
$$ M_ {U, \lambda } (x) =
\begin{cases}
\left ( \frac {\theta } {\theta + \lambda } \right ) ^ {S-x}  &   {\rm if} \ S-x \in \rm I\!N    \\
\left ( \frac {\theta } {\theta + \lambda } \right ) ^ {[S-x]+1}  &   {\rm if} \ S-x \notin \rm I\!N
\end{cases}
$$
where $[a]$ denotes the integer part of $a.$ Notice that
the condition $ M_ {U, \lambda } ( - \infty ) =0$ also holds.
Thus, recalling the expression of the Laplace transform of the Gamma density, we find that $\tau_S(x)$ has Gamma distribution
with parameters $(S-x, \theta)$ if $S-x$ is a positive integer, while it has Gamma distribution
with parameters $([S-x]+1, \theta)$ if $S-x$ is not an integer. This fact also follows directly, if one considers the nature of the Poisson process.
The moments
$ T_ n (x) = E \left [  \tau _S ^n (x) \right ]$ are soon obtained from the density or also by the formula  $T_n (x) =   (-1)^n \left [ \frac {\partial
^n } {\partial \lambda ^n } M_ {U, \lambda} (x) \right ] _ { \lambda =0}.$
We have:
$$ E ( \tau _S(x))=
\begin{cases}
\frac {S- x } {\theta } &  {\rm if} \ S- x \in \rm I\!N    \\
\frac {[S-x]+1 } {\theta}  &  {\rm if} \ S- x \notin \rm I\!N
\end{cases}
$$
and
$$ E ( \tau _S ^2(x))=
\begin{cases}
\frac {(S-x)^2 } {\theta ^2 } + \frac {S-x } {\theta ^2 }   &  {\rm if} \ S- x \in \rm I\!N    \\
\frac {([S-x]+1)^2 } {\theta ^ 2} + \frac {[S-x]+1 } {\theta ^2  }   &  {\rm if} \ S-x \notin \rm I\!N
\end{cases}
$$
Thus:
$$ Var( \tau _S(x))=
\begin{cases}
\frac {S-x } {\theta ^2 }  &  {\rm if} \ S-x \in \rm I\!N    \\
\frac {[S-x]+1 } {\theta ^2 }   &  {\rm if} \ S-x \notin \rm I\!N
\end{cases}
$$
\indent By Theorem \ref{laplacetransform} with $U(x)=x,$ we get the
Laplace transform $M_ {U, \lambda } (x)$ of $A_S(x)$ as the solution of the equation $L M _ {U,  \lambda } (x) = \lambda x M _ {U,  \lambda } (x),$ with
outer condition $M_ {U, \lambda } ( y) =1$ for $ y \ge S.$ By solving this equation, we get:
$$ M_ {U, \lambda } (x) =
\begin{cases}
\theta ^{S-x} \cdot \{(\theta + \lambda x) ( \theta +  \lambda (x+1)) \cdots (\theta +  \lambda (x+(S-x)-1))\}  ^ {-1}
 &   {\rm if} \ S-x \in \rm I\!N    \\
\theta ^ {[S-x]+1} \cdot \{(\theta + \lambda x) ( \theta + \lambda (x+1)) \cdots (\theta + \lambda (x + [S-x])) \}^ {-1} &   {\rm if} \ S-x \notin \rm I\!N
\end{cases}
$$
Notice that the condition $ M_ {U, \lambda } ( - \infty ) =0$ is
fulfilled. \par We observe that $M_ {U, \lambda } (x)$ turns out
to be the Laplace transform of a linear combination of $S-x $
independent exponential random variables with parameter $\theta,$
with coefficients $x, x+1, \dots , S-1,$ \  if $S-x$ is an
integer, while it is the Laplace transform of a linear combination
of $[S-x]+1 $ independent exponential random variables with
parameter $\theta,$ with coefficients $x, x+1, \dots , x+[S-x], $
if $S-x$ is not an integer. \par\noindent The $n-$th order moment
of $A_S(x)$ is given by $(-1)^n \left [ \frac {\partial ^n }
{\partial \lambda ^n } M_ {U, \lambda } (x) \right ] _ { \lambda
=0};$ calculating the first and second  derivative, after some
tedious computations,  we obtain
$$ E(A_S(x)) =
\begin{cases}
\frac {S-x } { 2 \theta} (x+S-1) & {\rm if } \ S-x =k  \in \rm I\!N   \\
\frac {[S-x]+1 } { 2 \theta} (2x+[S-x]) & {\rm if } \ S-x \notin \rm I\!N
\end{cases}
$$
and
$$ E ( A_S^2(x))=
\begin{cases}
\frac {k } { 12 \theta ^2} \{12x^2(k+1)+12x(k^2-1)+3k^3-2k^2-3k+2\} \ \ \ \ {\rm if } \  S-x =k \in \rm I\!N   \\
\frac {[S-x]+1 } { 12 \theta ^2} \{12x([S-x] +2) (x + [S-x]) +  \\ + [S-x] (3 [S-x] ^2 + 7 [S-x] +2 ) \}
  \ \ \ \ \ \ \ \ \ \ \ \ \ \ \ \ \ \ \ \ \ \ \ \ \ \ \ \ \ \ \ \ \  {\rm if } \ S-x \notin \rm I\!N
\end{cases}
$$
\begin{Remark}
{ \rm Notice that
$$ \tau _S(x)= \inf \{ t >0: x + N_t \ge S \} = \inf \{ t >0 : S -x - N_t \le 0 \} = \widehat \tau (S -x) $$
where $\widehat \tau (y)$ is the first-passage time below zero of the process $Y(t):= y - N_t .$ Moreover, it holds
$ A_S(x)= S \tau _S (x) - \widehat A (S-x),$ where $\widehat A (y)$ is the area swept out by $Y(t)$ till its first passage
below zero. Thus, the moments of
$\tau _S$ and the mean of $A _S $ can be obtained by the moments of $ \widehat \tau $ and $ E(\widehat A) ,$
which were calculated in \cite{abundo:passarea11}. In this way, one can avoid the heavy computations above.
}
\end{Remark}
\bigskip

\noindent{\bf Example 5} (a Levy process) \par
\noindent Let us consider the process $X(t)= x + \beta t + B_t + N_t ,$ where $N_t$ is a homogeneous Poisson Process with intensity $\theta >0$ and
let $S >0, \  \beta + \theta > S-x >0.$
We have $\tau _S(x)= \inf \{ t >0: X(t) \ge S \}  = \inf \{ t >0 : (S-x) - \beta t + W_t \le N_t \}$ i.e. the
first hitting time of BM $W_t$ with drift $ - \beta$ starting from $S- x,$ and the Poisson process $N_t$
(see  \cite{abundo:mcap08} as regards the density of $\tau _S(x))$ in a similar case). The condition
$ \beta + \theta > S-x $ assures that $\tau _S(x)$ is finite with probability one (see  Remark \ref{remark2.2}).
Now, the infinitesimal generator of the process is $Lf (x) = \frac 1 2 f''(x) + \beta f'(x) + f(x+1) - f(x) $ and
the differential-difference equations involved to find the Laplace transforms of $\tau _S(x)$ and $A_S(x),$ as well as those for
the moments of $\tau_S(x)$ and $A_S(x),$ cannot be solved explicitly; these quantities have to be found by a numerical procedure.

\begin{Remark}
{\rm
Until now we have supposed that the starting point $x < S$ is given and fixed.  One could introduce a randomness in the starting point, replacing $X(0)=x$ with a random variable $\eta ,$ having  density $g(x)$ whose support is the interval $ ( - \infty, S ).$ Thus, the quantities of interest become:
$ \tau _S = \inf \{ t >0: X (t) > S \} $ and
$ A_S = \int _ 0 ^{ \tau _S} X (t) dt ,$
while $\tau _S (x)$ and $A_S(x)$ are the corresponding values conditional to $ \eta = x .$ By using the found expressions for the moments of $\tau_S (x)$ and
$A_S(x),$ one can easily calculate the moments of $ \tau_S$ and $A_S.$  For instance, when $X(t)$ is  BM with drift $\mu,$ it follows by \eqref{momentstauforBM} and
\eqref{expressionmeanarea} that:
$$ E( \tau _S ^\mu )= \int _ {- \infty } ^S \frac { S - x } \mu \ g(x) dx = \frac 1 \mu E( S - \eta) $$
$$ E( (\tau _S ^\mu)^2 )= \int _ {- \infty } ^S  \left ( \frac { S - x } {\mu ^3} + \frac {(S-x)^2 } { \mu ^2 } \right ) \ g(x) dx =
\frac 1 {\mu ^3 } E( S - \eta) + \frac 1 { \mu ^2 } E((S- \eta ) ^2 )  $$
$$ E( A_S ^\mu ) = \int _ {- \infty } ^S \frac { S - x } { 2 \mu} \left ( S + x - \frac 1 \mu \right )  \ g(x) dx =
\frac 1 { 2 \mu } E \left [ \left ( S - \eta \right ) \left( S + \eta - \frac 1 \mu \right ) \right ]. $$
}
\end{Remark}

\section{Conclusions and Final Remarks}
For $S>0,$ we have considered a one-dimensional jump-diffusion process $X(t)$ starting from $x<S ,$  that is, a diffusion to which jumps
at Poisson-distributed instants are superimposed;
then, we have studied the probability distribution of the (random) area $A _S(x) = \int _0 ^{ \tau_S(x)} X(t) dt$ swept out by $X(t)$ till its
first-passage time  over the barrier $S,$ i.e. the
(random) time $\tau_S(x) = \inf \{ t >0 : X(t) \ge S | X(0)=x \} .$
The analogous problem concerning the first passage of $X(t)$ below zero was studied in \cite{abundo:passarea11}, while
results for Brownian motion with negative drift were obtained e.g. in
 \cite{janson:pro07}, \cite{kearney:jph05}, \cite{kearney:jph07},
\cite{knight:jam00}, \cite{perman:aap96}. \par\noindent
In particular, we have shown that the Laplace transforms of $A_S(x)$ and $\tau _S(x),$ their moments, as well as
the probability distribution of the
minimum  of $X(t)$ in $[0, \tau _S (x)],$ are solutions to certain partial differential-difference equations (PDDEs)
with outer conditions. Notice that, in the
absence of jumps, these PDDEs with outer conditions become simply PDEs with boundary conditions.
The quantities here investigated have interesting applications in Queueing Theory and in Economics (see the Introduction for a brief discussion). \par
After considering theoretical results for diffusions $X(t)$ with and without jumps, in the final part of the paper we  reported some
examples for which we have carried out explicit calculations. We remark that, in general it is not possible to solve explicitly the equations involved,
in order to find closed formulae for the Laplace transform of the first-crossing area, and for its moments;
moreover, even if one is able to find explicitly the Laplace transform, it is not always possible to invert it, to get the distribution of $A_S(x).$
\par
%The aim of this paper was to illustrate, also by means of
%practical examples, some techniques and results which  can be
%useful to study  the first-passage area of not only
%jump-diffusions, but also of more complex processes which are well
%approximated by  jump-diffusions, including  non-diffusion
%processes such as the piecewise-deterministic Markov pro-- \par
%\noindent cesses.
When the analytical solution is not available, due to the
complexity of calculations, one can resort to numerical solution
of the PDDEs involved; alternatively, one can carry out computer
simulation of a large enough number of trajectories of the process
$X(t),$ in order to obtain statistical estimations of the
quantities of interest.


\bigskip



\begin{thebibliography}{spc}

\bibitem [1] {abramo:handbook65}
Abramowitz, M.; Stegun, I. A. \newblock {\it Handbook of Mathematical Functions with Formulas, Graphs, and Mathematical Tables}.
\newblock Dover, New York, 1965

\bibitem [2] {abundo:passarea11}
Abundo, M. \newblock On the first-passage area of one-dimensional
jump-diffusion process.
\newblock{Methodol. Comput. Appl. Probab. Online First} 2011, DOI: 10.1007/s11009-011-9223-1.

\bibitem [3] {abundo:browarea11}
Abundo, Marco; Abundo, Mario \newblock
On the first-passage area of an emptying Brownian queue.
\newblock{International Journal of Applied Mathematics (IJAM)}
{\bf 2011}, {\it 24 (2)}, 259--266.


\bibitem [4] {abundo:asy09}
Abundo, M.\newblock
First-Passage Problems for Asymmetric
Diffusions and Skew-diffusion Processes.
\newblock{Open Systems $\&$ Information Dynamics} {\bf 2009}, {\it 16 (4)}, 325--350.

\bibitem [5] {abundo:ija08}
Abundo, M.\newblock
On the distribution of the time average of a jump-diffusion process.
\newblock{International Journal of Applied Mathematics} {\bf 2008}, {\it 21 (3)}, 447--454.

\bibitem[6]{abundo:euc07}
Abundo, M.\newblock
On first-passage problems for asymmetric one-dimensional diffusions. \newblock
{Lecture Notes in Computer Science, Computer Aided Systems Theory - EUROCAST 2007}; 4739, 179--186.\newblock
Springer Berlin/ Heidelberg, 2007.

\bibitem [7] {abundo:saa06}
Abundo, M.\newblock
Limit at zero of the first-passage time density and the inverse problem for
one-dimensional diffusions.
\newblock{Stochastic Anal. Appl.} {\bf 2006}, {\it 24}, 1119--1145.


\bibitem [8] {abundo:pms00}
Abundo, M.\newblock
On first-passage-times for one-dimensional jump-diffusion processes.
\newblock{Prob. Math. Statis.} {\bf 2000}, {\it 20 (2)}, 399--423.

\bibitem[9]{abundo:pms97}
Abundo, M.\newblock
On some properties of one-dimensional diffusion processes on an interval.
\newblock   {Prob.  Math.  Statis.} {\bf 1997}, {\it 17 (2)}, 235--268.


\bibitem[10]{abundo:mcap08}
Abundo, M.\newblock
On the First Hitting Time of a One-dimensional Diffusion and a Compound Poisson Process.
\newblock   {Methodol. Comput. Appl.  Probab.} {\bf 2010}, {\it 12}, 473--490

\bibitem[11]{darling:ams53}
Darling, D. A.; Siegert, A.J.F. \newblock
The first passage problem for a continuous Markov process.
\newblock   {Ann. Math. Statistics} {\bf 1953}, {\it 24}, 624--639.

\bibitem[12]
{gimsko:sde72}
Gihman, I.I.; Skorohod, A.V.\newblock
{\it Stochastic differential
equations.} \newblock
Springer-Verlag, Berlin, 1972.

\bibitem[13]
{grand:tab80}
Grandshteyn, I.S.; Ryzhik, I. M. \newblock
{\it Tables of Integrals, Series and Products.} \newblock
5th ed. Academic, London, 1980.

\bibitem[14]
{has:sto80}
Has' minskij, R.Z.\newblock
{\it Stochastic stability of differential
equations.} \newblock
Alphen a/d Rijn, Sijthoff $\&$ Noordhoff, 1980

\bibitem[15]
{ikwa:sde81}
Ikeda, N.; Watanabe, S.\newblock
{\it Stochastic differential equations and
diffusion processes.} \newblock
North-Holland Publishing Company, 1981

\bibitem[16]
{janson:pro07}
Janson, S.\newblock
Brownian excursion area, Wright's
constants in graph enumeration, and
other Brownian areas. \newblock {Probability Surveys} {\bf 2007}, {\it 4}, 80--145.

\bibitem[17]
{kartay:sto75}
Karlin, S.; Taylor, H.M.\newblock
{\it A second course in stochastic processes.} \newblock
Academic Press, New York, 1975.

\bibitem[18]
{kearney:jph05}
Kearney, M. J.; Majumdar, S.N.\newblock
On the area under a continuous time Brownian motion till its first-passage time. \newblock {J. Phys. A: Math. Gen.} {\bf 2005},
{\it 38}; 4097--4104.

\bibitem[19]
{kearney:jph07}
Kearney, M. J.; Majumdar, S.N.; Martin R.J.\newblock
The first-passage area for drifted Brownian motion and the moments of the Airy distribution.
\newblock {J. Phys. A: Math. Theor.} {\bf 2007}, {\it 40}, F863--F864.

\bibitem [20] {klebaner:99}
Klebaner, F. C.\newblock
\ {\it Introduction to stochastic calculus with applications.}
\newblock Imperial College Press, Singapore, 1999.


\bibitem[21]
{knight:jam00}
Knight, F. B. \newblock
The moments of the area under reflected Brownian Bridge conditional on its local time at zero.
\newblock {Journal of Applied Mathematics and Stochastic Analysis} {\bf 2000},
{\it 13 (2)}, 99--124.

\bibitem[22]{lansky:jtb94}
Lanska, V.; Lansky, P.; Smiths, C.E. \newblock
Synaptic transmission in a diffusion model for neural activity.
\newblock   {J. Theor. Biol. } {\bf 1994}, {\it 166}, 393--406.


\bibitem[23]
{perman:aap96}
Perman, M.; Wellner, J. A. \newblock
On the Distribution of Brownian Areas. \newblock
{The Annals of Applied Probability} {\bf 1996}, {\it 6 (4)}, 1091--1111.

\bibitem[24]
{ricc:jap99}
Ricciardi, L.M; Di Crescenzo, A.; Giorno V.; Nobile, A.G. \newblock
An outline of theoretical and algorithmic approaches to first passage time problems with applications to biological modeling. \newblock
{Math. Japonica} {\bf 1999}, {\it 50 (2)}, 247--322.

\bibitem[25]
{tuckwell:jap76}
Tuckwell, H.C. \newblock
On the first exit time problem for temporally homogeneous Markov processes. \newblock
{The Annals of Applied Probability} {\bf 1976}, {\it 13}, 39--48.


\end{thebibliography}
\end{document}